\documentclass[12pt] {article}
\usepackage{amsmath,amsthm}
\usepackage{amssymb,latexsym}
\usepackage{amscd}

\title{On repeated sequential closures of constructible functions in valuations.}
\date{}
\author{ Semyon Alesker \footnote{Partially supported by ISF grant 1447/12.}
\\  { \normalsize Department of Mathematics, Tel Aviv University, Ramat Aviv}
 \\  { \normalsize 69978 Tel Aviv,
Israel }
\\ {\normalsize e-mail: semyon@post.tau.ac.il}}

\def\eps{\varepsilon}
\def\alp{\alpha}
\def\ome{\omega}

\def\lam{\lambda}

\def\to{\rightarrow}
\def\qed {Q.E.D.}
\def\pt{\partial}

\def\RR{\mathbb{R}}
\def\CC{\mathbb{C}}

\def\NN{\mathbb{N}}
\def\ZZ{\mathbb{Z}}

\def\PP{\mathbb{P}}

\def\One{{1\hskip-2.5pt{\rm l}}}
\swapnumbers
\newtheorem{theorem}{Theorem}[section]
\newtheorem{corollary}[theorem]{Corollary}
\newtheorem{lemma}[theorem]{Lemma}
\newtheorem{proposition}[theorem]{Proposition}
\newtheorem{claim}[theorem]{Claim}
\theoremstyle{definition}

\theoremstyle{proposition-definition}
\newtheorem{proposition-definition}[theorem]{Proposition-Definition}

\numberwithin{equation}{section}

  \def\cf{{\cal F}}
  
 \def\ck{{\cal K}}

  \def\cx{{\cal X}}

\newcommand \supp{\operatorname{supp} \,}

\begin{document}
\def\bgl{\overline{GL}}
\def\secF{\tilde\cf}

\maketitle

\begin{abstract}
The space of constructible functions form a dense subspace of the space of generalized valuations. In this note we prove a somewhat stronger property that
the sequential closure, taken sufficiently many (in fact, infinitely many) times, of the former space is equal to the latter one. This stronger property is necessary for some applications
in \cite{alesker-valuations-and-integral-geometry}.
\end{abstract}

\section{Main results.}\label{S:main-results}
The main results of this note are Theorem \ref{T:1} and Corollaries \ref{C:1-1}, \ref{C:1-2} below. Corollary \ref{C:1-1} says that the taken sufficiently (infinitely) many times
sequential closure of constructible functions inside the space of generalized valuations is equal to the whole space.
Corollary \ref{C:1-2} says, in particular, that if a sequentially continuous
linear operator from generalized valuations on a manifold with or without compact support to a Hausdorff linear topological space vanishes on constructible functions, then it vanishes.
Recall that a map between two topological spaces is called sequentially continuous if it maps convergent sequences to convergent ones. Notice that for non-metrizable topological spaces sequential continuity of a map
does not imply topological continuity.

The reason to write this note is to correct a mistake made by the author in \cite{alesker-valuations-and-integral-geometry}, where it was wrongly claimed that
several operations, such as pull-back, push-forward, and product on generalized valuations with given wave front sets are topologically continuous, while they satisfy, in fact,
only a weaker property of sequential continuity. This property comes from the fact that operations of pull-back, push-forward, and product on generalized functions
or distributions with given wave front sets are only sequentially continuous in appropriate (so called H\"ormander) topology, but not topologically continuous.\footnote{This fact was pointed out to the author
by C. Brouder in September 2013. I am very grateful to him for this remark.}

\hfill

Let $X$ be either smooth manifold or real analytic manifold. It will always be assumed to be countable at infinity,
i.e. $X$ can be presented as a countable union of compact subsets. We will denote by $V^\infty_c(X)$ and $V^\infty(X)$ the space of smooth valuations on $X$ with and
without compact support respectively.
Also we denote by $V^{-\infty}_c(X)$ and $V^{-\infty}(X)$ the spaces of generalized valuations with and without compact support respectively.
We refer to \cite{alesker-valuations-and-integral-geometry} for the definitions of all these spaces and their topologies (see also \cite{alesker-gafa-survey}).

Furthermore, by the abuse of notation, for a smooth manifold $X$ we denote by $\cf(X)$ the
linear span of indicator functions of compact submanifolds with corners; while for a real analytic $X$ we denote by $\cf(X)$ the space of so called $\CC$-valued constuctible functions.
Let us remind the definition of this notion following \cite{part4}.

We refer to \S 8.2 in \cite{kashiwara-schapira} for the definition and basic properties of subanalytic sets (see also Section 1.2 in \cite{part4}).
An integer valued function $f\colon X\to \ZZ$ on a real analytic manifold is called constructible if it satisfies:

1) for every $m\in \ZZ$ the set $f^{-1}(m)$ is subanalytic;

2) the family of sets $\{f^{-1}(m)\}_{m\in \ZZ}$ is locally finite.

Now a $\CC$-valued function $f\colon X\to \CC$ is called {\itshape constructible} if $f$ is a finite linear combination
with $\CC$-coefficients of integer valued constructible (in the above sense) functions.

\hfill

Furthermore, for a closed subset $Z\subset X$ we denote by $\cf_Z(X)$ one of the following:

(1) if $X$ is a smooth manifold, then $\cf_Z(X)$ is the $\CC$-linear span of
indicator functions of compact submanifolds $P\subset X$ with corners  such that $P\subset Z$;

2) if $X$ is a real analytic manifold, then $\cf_Z(X)$ is the subspace of $\CC$-valued constructible functions
supported in $Z$.

\hfill

In order to formulate the main results let us remind the notion of sequential closure of transfinite order of a set.
This notion was already known to S. Banach, see p. 213 in his classical book \cite{banach}.
Let $\cx$ be a topological space. Let $A\subset \cx$ be a subset. A sequential closure of $A$ is defined by
$$scl(A)=\{x\in \cx|\, \exists \mbox{ sequence } \{a_i\}\subset A,\, a_i\to x\}.$$
It is clear that $A\subset scl(A)$, and if $\cx$ is a linear topological space and $A\subset \cx$ is a linear subspace
then $scl(A)$ is a linear subspace. In general, if $\cx$ is not metrizable, $scl(A)$ may not be closed or even sequentially closed, i.e.
$scl(scl(A))\ne scl(A)$. We can repeat the procedure of taking sequential closure any number of times, even any infinite number of times corresponding to any ordinal.
More precisely, for any ordinal $\eta$ one can define by transfinite induction the subset $scl^\eta(A)$ as follows:

$\bullet$ if $\eta=0$ then $scl^0(A)=A$;

$\bullet$ if $\eta=\xi+1$ then $scl^\eta(A)=scl(scl^\xi(A))$;

$\bullet$ if $\eta$ is a limit ordinal then $scl^\eta(A)=\cup_{\xi<\eta}scl^\xi(A)$.

Furthermore there exists an ordinal $\eta$ such that for any $\eta'>\eta$ one has $scl^{\eta'}(A)=scl^\eta(A)$. We will denote the latter subset by $scl^*(A)$.
It is also clear that $scl^*(A)$ is sequentially closed, i.e. $scl(scl^*(A))=scl^*(A)$.
Clearly if $\cx$ is a linear topological space and $A\subset \cx$ is a linear subspace then $scl^\eta(A)$ is a linear subspace for any $\eta$.

\hfill

Here is the first main result, where for a compact subset $Z\subset X$ we denote by $V^{-\infty}_Z(X)$ the space of generalized valuations
with support contained in $Z$.

\begin{theorem}\label{T:1}
Let $X$ be either real analytic or smooth manifold countable at infinity. Let $Z_1\subset X$ be a compact
subset. Let $Z_2$ be a compact neighborhood of $Z_1$. Then in the above notation the subspace
$scl^*(\cf_{Z_2}(X))\subset V_c^{-\infty}(X)$ contains $V^{-\infty}_{Z_1}(X)$.
\end{theorem}

Let us deduce two immediate corollaries.
\begin{corollary}\label{C:1-1}
Let $X$ be either real analytic or smooth manifold countable at infinity. Then in the above notation
$scl^*(\cf(X))=V^{-\infty}(X)$.
\end{corollary}
{\bf Proof of Corollary \ref{C:1-1}.} By Theorem \ref{T:1}
$scl^*(\cf(X))$ contains $V^{-\infty}_c(X)$. But $V^{-\infty}_c(X)$ is
sequentially dense in $V^{-\infty}(X)$ since $X$ is assumed to be
countable at infinity. \qed

\begin{corollary}\label{C:1-2}
Let $X$ be either real analytic or smooth manifold countable at infinity. Let
$R\colon V^{-\infty}(X)(\mbox{resp. }V_c^{-\infty}(X))\to E$ be a linear operator into a Hausdorff topological vector space $E$. Assume that $R$
is sequentially continuous, i.e. $R$ maps convergent sequences in $V^{-\infty}(X)(\mbox{resp. }V^{-\infty}_c(X))$ to convergent sequences in $E$.
Let $L\subset E$ be a sequentially closed subset.\footnote{$L\subset E$ is called sequentially closed if for any converging in $E$ sequence $\{x_i\}_{i=1}^\infty\subset L$ its limit
belongs to $L$; equivalently $scl(L)=L$. Clearly any closed subset is sequentially closed, but the converse is not true in general.}
Then if $R(\cf(X))\subset L$ then the whole image of $R$ is contained in $L$. In particular if $R(\cf(X))=0$ then $R\equiv 0$.
\end{corollary}
{\bf Proof.} By transfinite induction $R(scl^*(\cf(X)))\subset L$. In all cases (smooth or real analytic $X$, compact or non-compact support) $scl^*(\cf(X))$
is equal to the whole space. \qed

\section{Proofs.}\label{S:proofs}
The two cases of real analytic and smooth manifold $X$ are almost identical and will be treated simultaneously.

\begin{lemma}\label{L:smooth-val}
Let $\sigma$ be an infinitely smooth measure on a vector space $V$. Let $A\in \ck^\infty(V)$.
Then $K\mapsto \sigma(K+A)$ is a smooth valuation.
\end{lemma}
{\bf Proof.} Consider the map
$$p\colon V\times \PP_+(V^*)\times[0,1]\to V$$
given by $(x,n,t)=x+t\nabla h_A(n)$, where $h_A\colon V^*\to \RR$ is the supporting functional of $A$.
Since $h_A$ is 1-homogeneous, its gradient $\nabla h_A$ is 0-homogeneous, and hence considered as a map
$\nabla h_A\colon \PP_+(V^*)\to V$. Since $A\in \ck^\infty(V)$ the latter map is infinitely smooth.

We may and will assume that $0\in int(A)$; the general case reduces to this one by a translation. In this case
the restriction of $p$ to $N(K)\times [0,1]$ is a homeomorphism onto the closure of $(K+A)\backslash K$. Hence
\begin{eqnarray}\label{E:j-1}
\sigma(K+A)=\sigma(K)+\int_{N(K)\times[0,1]}p^*\sigma.
\end{eqnarray}
Let $q\colon V\times \PP_+(V^*)\times[0,1]\to V\times \PP_+(V^*)$ is the obvious projection.
Then by (\ref{E:j-1}) we have
$$\sigma(K+A)=\sigma(K)+\int_{N(K)}q_*p^*\sigma.$$
Obviously $q_*p^*\sigma$ is a smooth $(\dim V-1)$-form on $V\times \PP_+(V^*)$. This proves the lemma.


\hfill

\begin{lemma}\label{L:01}
Let $Z_1\subset V$ be a compact set, and $Z_2$ be a compact
neighborhood of $Z_1$. Let $\sigma$ be a smooth measure with
$supp(\sigma)\subset Z_1$. Let $A\in \ck^\infty_{an}(V)$ such that
$Z_1-A\subset Z_2$. Then

(1) one has
\begin{eqnarray}\label{E:01}
\sigma(\bullet +A)=\int \One_{x-A}\cdot d\sigma(x)
\end{eqnarray}
as generalized valuations;

(2) $\sigma(\bullet +A)$ belongs to $V^\infty_{Z_2}(V)\cap
scl^*(\cf_{Z_2}(V))$.
\end{lemma}
{\bf Proof.} First notice that $\sigma(\bullet +A)$ is a smooth valuation by Lemma \ref{L:smooth-val}.

To prove part (1) it suffices to apply both sides to an arbitrary
smooth compactly supported valuation and to prove that the result is
the same. It suffices to apply them to such valuations of the form
$\eta=\ome(\bullet+B)$ where $\ome$ is a smooth compactly supported
measure, and $B\in \ck^\infty(V)$.

Apply the left hand side of (\ref{E:01}) to $\eta$:
\begin{eqnarray*}\label{E:01.5}
<\sigma(\bullet +A),\eta>=(\sigma\boxtimes\ome)(\Delta(V)+(A\times
B))\overset{Fubini}{=}\\
\int d\sigma(x)\ome(x-A+B).
\end{eqnarray*}
Apply the right hand side of (\ref{E:01}) to $\eta$:
\begin{eqnarray*}
<\int d\sigma(x) \One_{x-A},\eta>=\int d\sigma\cdot \eta(x-A)=\int
d\sigma\cdot \ome(x-A+B)=\\<\sigma(\bullet +A),\eta>.
\end{eqnarray*}
Thus part (1) is proved.

It remains to show that $\sigma(\bullet
+A)\in scl^*(\cf_{Z_2}(V))$. We use the equality (\ref{E:01}) and replace
the integral in the right hand side by a Riemann sum corresponding
to a subdivision of $V$ whose diameter will tend to 0. Let us show
that these Riemann sums converge to the integral in the weak
topology on $V^{-\infty}(V)$. Let $\{C_i^N\}_{N=1}^\infty$ be a
sequence of subdivisions whose diameter tends to 0 as $N\to \infty$.
Choose a point $x_i^N\in C_i^N$. Apply the corresponding Riemann sum
to $\eta$:
\begin{eqnarray*}
<\sum_i
\sigma(C_i^N)\One_{x_i^N-A},\eta>=\\
\sum_i\sigma(C_i^N)\eta(x_i^N-A)\underset{N\to \infty}{\to}\int
d\sigma(x)\eta(x-A)=\\
<\int d\sigma(x) \One_{x-A}, \eta>,
\end{eqnarray*}
where we have applied the fact that continuous scalar valued
functions are Riemann integrable,  to the function $[x\mapsto
\eta(x-A)]$. It only remains to notice that $i$-th Riemann sum
belongs to $\cf_{Z_2}(V)$ for $i\gg 1$.\qed

\begin{lemma}\label{L:02}
Let $Z_1\subset V$ be a compact set, and $Z_2$ be a compact
neighborhood of $Z_1$. Let $A_1,A_2\in \ck^\infty_{an}(V)$ such that
$Z_1-A_1\subset Z_2$. Let $\sigma_1,\sigma_2$ be smooth measures on
$V$ such that $supp(\sigma_1)\subset Z_1$. Let
$\phi_i:=\sigma_i(\bullet+A_i)$, $i=1,2$. Then
$$\phi_1\cdot\phi_2\in V^\infty_{Z_2}(V)\cap scl^*(\cf_{Z_2}(V)).$$
\end{lemma}
{\bf Proof.} First let us show that in the space $V^{-\infty}(V)$
\begin{eqnarray}\label{E:02}
\phi_1\cdot\phi_2=\int\int d\sigma_1(x)
d\sigma_2(y)\One_{(x-A_1)\cap(y-A_2)},
\end{eqnarray}
where the integral is understood in the sense the limit of Riemann sums converging in the weak topology.
Again we have to show that if we apply the two sides on the same
$\eta\in V^\infty_c(V)$ then we get the same result. It suffices to
choose $\eta$ of the form $\eta=\ome(\bullet +C)$, where
$C\in\ck^\infty(V)$, and $\ome$ is a smooth compactly supported
measure on $V$. Applying the right hand side of (\ref{E:02}) on such
$\eta$ we get
\begin{eqnarray}
<\int\int d\sigma_1(x)d\sigma_2(y)\One_{(x-A_1)\cap(y-A_2)},\eta>=\\
\int\int
d\sigma_1(x)d\sigma_2(y)\eta((x-A_1)\cap(y-A_2))=\\\label{E:03}
\int\int
d\sigma_1(x)d\sigma_2(y)\ome\left([(x-A_1)\cap(y-A_2)]+C\right).
\end{eqnarray}

Now let us apply to $\eta$ the left hand side of (\ref{E:02}) (in
the computation $\Delta$ is the diagonal map $V\to V\times V\times
V$ given by $x\mapsto (x,x,x)$):
\begin{eqnarray}
<\phi_1\cdot\phi_2,\eta>=(\phi_1\cdot\phi_2\cdot\eta)(V)=\\
(\sigma_1\boxtimes\sigma_2\boxtimes\ome)(\Delta(V)+(A_1\times
A_2\times C))\overset{Fubini}{=}\\
\int\int d\sigma_1(x)d\sigma_2(y)\ome([(x-A_1)\cap
(y-A_2)]+C)=(\ref{E:03}).
\end{eqnarray}
Thus equality (\ref{E:02}) is proven. It remains to show that the
right hand side of (\ref{E:02}) belongs to $scl^*(\cf_{Z_2}(V))$. To do
that, we will approximate the double integral by Riemann sums
belonging to $\cf_{Z_2}(V)$ which converge to the double integral in
the weak topology on $V^{-\infty}(V)$.

Consider a sequence $\{C_i^N\}_{N=1}^\infty$ of subdivisions of $V$
with diameter tending to 0 as $N\to \infty$. Choose a point
$x_i^N\in C_i^N$. For the corresponding Riemann sum
\begin{eqnarray}
<\sum_{i,j}\sigma_1(C_i^N)\sigma_2(C_j^N)\One_{(x_i^N-A_1)\cap
(x_j^N-A_2)},\eta>=\\\label{E:double-riemann}
\sum_{i,j}\sigma_1(C_i^N)\sigma_2(C_j^N)\eta((x_i^N-A_1)\cap
(x_j^N-A_2)).
\end{eqnarray}

Similarly for the double integral we have
\begin{eqnarray}
<\int\int d\sigma_1(x)
d\sigma_2(y)\One_{(x-A_1)\cap(y-A_2)},\eta>=\\\label{E:double-int}
\int\int d\sigma_1(x) d\sigma_2(y)\cdot\eta((x-A_1)\cap(y-A_2)).
\end{eqnarray}
We see that (\ref{E:double-riemann}) is a Riemann sum for
(\ref{E:double-int}), and we have to show that the former converges
to the latter. In other words we have to show that the function
$V\times V\to\RR$ given by $(x,y)\mapsto \eta((x-A_1)\cap(y-A_2))$
is Riemann integrable.

Notice that the above function does not have to be continuous. But
obviously this function is bounded. By the Lebesgue criterion of
Riemann integrability (see e.g. \cite{zorich-book}, Section 11.1) it suffices to show that
the above function is continuous almost everywhere. For that it
suffices to prove that the function $\Xi\colon V\times V\to
\ck(V)\cup\{\emptyset\}$ given by $\Xi(x,y)=(x-A_1)\cap (y-A_2)$ is
continuous almost every where is the Hausdorff metric on $\ck(V)$.

To prove the last statement let us consider a closed convex set
$M\subset V\times V\times V$ defined by
$$M:=\{(x,y,z)|\, x-z\in A_1,y-z\in A_2\}.$$
Let $q\colon V\times V\times V\to V\times V$ be the projection onto
the first two copies of $V$. Then clearly
$$q^{-1}(x,y)\cap M=(x-A_1)\cap (y-A_2)=\Xi(x,y),$$
and the restriction of $q$ to $M$ is proper. Applying Theorem 1.8.8
of \cite{schneider-book}, it follows that $\Xi$ is continuous
outside of the boundary of the set
$$\Xi(M)=\{(x,y)|(x-A_1)\cap (y-A_2)\ne \emptyset\}=\{(x,y)| x-y\in
A_1-A_2\}.$$

From this description it is clear that $\Xi(M)$ is a closed convex
set. Its boundary always has Lebesgue measure zero. Finally let us
notice that $N$-th Riemann sum belongs to $\cf_{Z_2}(V)$ for $N\gg
1$. \qed

\begin{lemma}\label{L:001}
Let $Z_1\subset V$ be a compact domain with infinitely smooth
boundary, and $Z_2$ be a compact neighborhood of $Z_1$. The image of
$W_{i,Z_1}^\infty\cap scl^*(\cf_{Z_2}(V))$ in
$W^\infty_{i,Z_1}/W^\infty_{i+1,Z_1} \simeq
C_{Z_1}^\infty(V,Val_i^\infty(V))$\footnote{This canonical
isomorphism was proved in Lemma 5.1.3(1) of \cite{part4}.}
is a dense linear subspace.
\end{lemma}
{\bf Proof.} Let us fix $A\in \ck^\infty_{an}(V)$. Let
$\phi'_\eps(K):=vol(K+\eps A)$. Define
$$\phi(K):=\frac{i!}{n!}\frac{d^{n-i}}{d\eps^{n-i}}\big|_{\eps=0}\phi'_\eps(K)=V(K[i],A[n-i]).$$
Clearly $\phi\in W_i^\infty$. Now let $\psi=\ome(\bullet+B)$, where
$B\in\ck^\infty_{an}(V)$ such that $0\in B$, and $\ome$ is a smooth
compactly supported measure on $V$ such that $supp(\ome)+B\subset
Z_1$. Thus $supp(\psi)\subset Z_1$. By Lemma \ref{L:02} for small $\eps>0$
$$\psi\cdot\phi'_\eps\in V^\infty_{Z_1}(V)\cap scl^*(\cf_{Z_2}(V)).$$
Hence also $\psi\cdot \phi\in V^\infty_{Z_1}(V)\cap scl^*(\cf_{Z_2}(V))$.
But since $\phi\in W_i^\infty$ and $V^\infty(V)\cdot
W_i^\infty\subset W_i^\infty$, we deduce that
\begin{eqnarray}\label{E:1-1}
\phi\cdot \psi\in W_{i,Z_1}^\infty\cap scl^*(\cf_{Z_2}(V)).
\end{eqnarray}

Next let us compute the image of this valuation in
$C_{Z_1}^\infty(V,Val_i^\infty(V))$. We have
$$(\phi\cdot\psi)(K)=\frac{i!}{n!}\frac{d^{n-i}}{d\eps^{n-i}}|_{\eps=0}(\ome\boxtimes
vol)(\Delta(K)+(B\times\eps A)).$$

For any valuation $\xi\in W^\infty_i$, $K\in \ck(V), x\in V$
one has
$$\xi(x+\lam K)=O(\lam^i) \mbox{ as } \lam\to +0;$$
this was the definition of $W^\infty_i$ in \cite{part1}, beginning of Section 3 (where it was denoted by $W_i$).

For such a $\xi$, its image $\bar\xi$ in
$W_i^\infty/W^\infty_{i+1}=C^\infty(V,Val_i^\infty(V))$ is computed
as follows:
\begin{eqnarray*}
(\bar\xi(x))(K)=\underset{\lam\to+0}\lim \frac{1}{\lam^i}\xi(x+\lam
K).
\end{eqnarray*}
The limit necessarily exists and the map $\bar\xi$ takes values in
$Val_i^\infty(V)$.

For $\xi=\phi\cdot\psi$ as above we have
\begin{eqnarray}
\overline{\phi\cdot\psi}(x)(K)=\frac{1}{n!}\frac{\pt^n}{\pt\lam^{i}\pt\eps^{n-i}}\big|_{\lam=\eps=0}
(\ome\boxtimes vol)(\Delta(\lam K+x)+(B\times\eps A))=\\
\frac{1}{n!}\frac{\pt^n}{\pt\lam^{i}\pt\eps^{n-i}}\big|_{\lam=\eps=0}
(\ome\boxtimes vol)(\Delta(\lam K)+((x+B)\times(x+\eps A)))=\\\label{E:001}
\frac{1}{n!}\frac{\pt^n}{\pt\lam^{i}\pt\eps^{n-i}}\big|_{\lam=\eps=0}
((T_{-x})_*\ome\boxtimes vol)(\Delta(\lam K)+(B\times\eps A)),
\end{eqnarray}
where $(T_{-x})_*$ denotes the push-forward on measures under the
shift by $-x$, namely $[y\mapsto y-x]$. We will need a lemma.

\begin{lemma}\label{L:002}
Let $\sigma$ be a smooth measure on an $n$-dimensional vector space
$V$. Let $A,B,K\in \ck(V)$. Define the function of $(\lam,\eps)\in
[0,\infty)^2$ by
$$F_\sigma(\lam,\eps):=(\sigma\boxtimes vol)(\Delta(\lam K)+(B\times \eps
A)).$$ Then the following holds:

(1) $F_\sigma\in C^\infty([0,\infty)^2).$

(2) For $0\leq i\leq n$ define
$$h_{\sigma,i}(\lam):=\frac{i!}{n!}\frac{\pt^{n-i}}{\pt\eps^{n-i}}\big|_{\eps=0}F_\sigma(\lam,\eps).$$
Then $h_{\sigma,i}(\lam)=O(\lam^i)$ as $\lam\to +0$.

(3) $\underset{\lam\to
+0}{\lim}\frac{h_{\sigma,i}(\lam)}{\lam^i}=\sigma(B)\cdot
V(K[i],A[n-i])$.
\end{lemma}
Let us postpone the proof of Lemma \ref{L:002} and finish the proof
of Lemma \ref{L:001}. By Lemma \ref{L:002} and (\ref{E:001}) we have
\begin{eqnarray*}
\overline{\phi\cdot\psi}(x)(K)=\ome((T_{-x})_*B)\cdot
V(K[i],A[n-i])=\ome(x+B)\cdot V(K[i],A[n-i]).
\end{eqnarray*}
To summarize, we have proven so far the following: any smooth
$Val^\infty_i(V)$-valued function on $V$ of the form
$$x\mapsto \ome(x+B)\cdot V(\bullet[i],A[n-i]),$$
belongs to the image of $W_{i,Z_1}^\infty\cap scl^*(\cf_{Z_2}(V))$ in
$W^\infty_{i,Z_1}/W_{i+1,Z_1}^\infty$, where $A,B\in \ck^\infty(V)$
such that $0\in B$, and $\ome$ is a smooth measure on $V$ such that
$supp(\ome)+B\subset Z_1$. Now let us show that the closure of such
functions in the usual Fr\'echet topology on
$C^\infty_{Z_1}(V,Val_i^\infty(V))$ is equal to the whole space.

Let $B$ be the unit Euclidean ball in $V$. For any $l\in \NN$ the
function $\frac{\ome(x+\frac{1}{l}B)}{vol(\frac{1}{l}B)}\cdot
V(\bullet[i],A[n-i])$ belongs to the image of $W_i^\infty\cap
scl^*(\cf_{Z_2}(V))$. However obviously
$$\frac{\ome(x+\frac{1}{l}B)}{vol(\frac{1}{l}B)}\to
\frac{\ome}{vol}(x)\mbox{ in } C^\infty_{Z_1}(V).$$

This implies that for any smooth function $h\colon V\to \CC$, any
$A\in \ck^\infty(V)$ with $supp(h)\subset int(Z_1)$ the function
\begin{eqnarray}\label{E:002}
[x\mapsto h(x)\cdot V(\bullet[i],A[n-i])]
\end{eqnarray}
belongs to the closure of the image of $W_{i,Z_1}^\infty\cap
scl^*(\cf_{Z_2}(V))$. Since $i$-homogeneous mixed volumes are dense in
$Val_i^\infty(V)$ by \cite{alesker-mcmullen} we deduce that for any $h\in
C^\infty_{Z_1}(V)$ and any $\mu\in Val_i^\infty(V)$ the
$Val_i^\infty(V)$-valued function $h\otimes \mu$ lies in the closure
of the image of $W_{i,Z_1}^\infty\cap scl^*(\cf_{Z_2}(V))$ (here it is
the only place where we have used that the boundary of $Z_1$ is
smooth). But linear combinations of such elements are dense in
$C^\infty_{Z_1}(V,Val_i^\infty(V))$.  \qed

\begin{lemma}\label{L:approximation}
Let $X$ be a smooth manifold and $Z\subset X$ be a compact subset. Let $Z'\subset X$ be a
compact neighborhood of $Z$. Then for any element $\psi\in V^{-\infty}_Z(X)$ there exists a sequence
of elements from $V^{\infty}_{Z'}(X)$ converging to $\psi$ in the topology of $V^{-\infty}_c(X)$ (or equivalently
in the weak topology on $V^{-\infty}(X)$).
\end{lemma}
{\bf Proof.} \underline{Step 1.} Let us prove the statement for $X=\RR^n$. For let us choose a sequence $\{\mu_i\}$
of smooth non-negative compactly supported measures on the Lie group $Aff(\RR^n)$ of affine transformations of $\RR^n$ such that
$\int_{Aff(\RR^n)}\mu_i=1$ and $supp\{\mu_i\}\to \{id\}$ in Hausdorff metric on $Aff(\RR^n)$. Define
$$\psi_i:=\int_{g\in Aff(\RR^n)}g^*(\psi)\cdot d\mu_i(g).$$
In the proof of Lemma 8.2 in \cite{alesker-bernig} it was shown that $\psi_i\in V^\infty(X)$ for all $i$ and
$\psi_i\to \psi$ in $V^{-\infty}(X)$. It is also clear that for $i\gg 1$ one has $supp(\psi_i)\subset Z'$. This implies the
lemma for $X=\RR^n$.

\hfill

\underline{Step 2.} Assume now that $X$ is a general manifold. Let us choose a finite open covering $\{V_\alp\}_{\alp}$
of $Z$ and open subsets $U_\alp$ such that $V_\alp\subset U_\alp\subset Z'$, the closures $\bar U_\alp$ are compact and are contained in
the interior of $Z'$, and there exist diffeomorphisms $U_\alp\tilde\to \RR^n$.

Let us choose a partition of unity in valuations subordinate to the covering $\{V_\alp\}_\alp\cup\{X\backslash Z\}$; we denote it by $\{\xi_\alp\}_\alp\cup\{\xi\}$
where $supp(\xi_\alp)\subset V_\alp,\, \supp(\xi)\subset X\backslash Z$, and $\sum_\alp\xi_\alp+\xi=\chi$, where $\chi$ is the Euler characteristic.
Such partition of unity exists by \cite{part4}, Proposition 6.2.1. Since $supp(\psi)\subset Z$ we have
$$\psi=\sum_\alp\psi\cdot \xi_\alp.$$
Since $\psi\cdot\xi_\alp$ has compact support contained in $U_\alp$ and $U_\alp\simeq \RR^n$, Step 1 implies that that there exists
a sequence $\{\psi_{\alp,i}\}_i\subset V^\infty_{\bar U_\alp}(X)$ converging to $\psi\cdot \xi_\alp$ in the topology of $V^{-\infty}(X)$ as $i\to \infty$.
Then the sequence
$$\psi_i:=\sum_\alp\psi_{\alp,i}$$
satisfies the proposition. \qed

\hfill

{\bf Proof of the main result.} We have to show that
$V^{-\infty}_{Z_1}(X)\subset scl^*(\cf_{Z_2}(X))$ for a smooth manifold
$X$. Let us fix a compact neighborhood $Z'_1$ of $Z_1$ with
infinitely smooth boundary contained in the interior of $Z_2$. By Lemma \ref{L:approximation}
for every element of $V_{Z_1}^{-\infty}(X)$ there
exists a sequence of elements of $V^\infty_{Z'_1}(X)$ converging to
this element in the weak topology on $V^{-\infty}(X)$. Hence it suffices
to show that
\begin{eqnarray}\label{E:003}
V^\infty_{Z_1'}(X)\subset scl^*(\cf_{Z_2}(X)).
\end{eqnarray}
Notice that $V^\infty_{Z'_1}(X)\cap scl^*(\cf_{Z_2}(X))$ is a closed
subspace of $V^\infty_{Z'_1}(X)$ since $V^\infty_{Z'_1}(X)$ is
metrizable.

First let us prove (\ref{E:003}) for $X$ being a vector space. If
this is not true then there exists a unique integer $0\leq i\leq n$
such that $W^\infty_{i+1, Z'_1}\cap
scl^*(\cf_{Z_2}(X))=W^\infty_{i+1,Z'_1}$ and $W^\infty_{i,Z'_1}\cap
scl^*(\cf_{Z_2}(X))\ne W^\infty_{i,Z'_1}$. In this case the image of
$W^\infty_{i,Z'_1}\cap scl^*(\cf_{Z_2}(X))$ in
$W_{i,Z'_1}^\infty/W_{i+1,Z'_1}^\infty$ is a closed subspace.
However by Lemma \ref{L:001} this image is dense. Hence
$W^\infty_{i,Z'_1}\cap scl^*(\cf_{Z_2}(X))= W^\infty_{i,Z'_1}$ which is a
contradiction. This proves (\ref{E:003}) when $X$ is a vector space.

Let us prove (\ref{E:003}) for a general manifold $X$. Let us fix a
finite open covering $\{U_\alp\}$ of $Z'_1$ such that the closures
$\bar U_\alp\subset int(Z_2)$ and each $U_\alp$ is analytically
diffeomorphic to $\RR^n$. By Proposition 6.2.1 of
\cite{part4} one can construct a partition of unity in
valuations subordinate to this covering, namely there exist
valuations $\{\phi_\alpha\}$ such that $supp(\phi_\alpha)\subset
U_\alp$ and in a neighborhood of $Z_1'$ one has $\sum_\alp\phi_\alp=\chi$ (here $\chi$ is the Euler
characteristic). Any $\psi\in V^\infty_{Z'_1}(X)$ can be written
$\psi=\sum_\alp\phi_\alp\cdot \psi$. Let us choose compact sets
$Z_{\alp,2}\subset U_\alp$ such that $supp(\phi_\alp)\subset
int(Z_{\alp,2})$ and $Z_{\alp,2}\subset Z_2$. As we have shown for a
vector space,
$$\phi_\alp\cdot \psi\in scl^*(\cf_{Z_{\alp,2}}(U_\alp)).$$
But obviously $scl^*(\cf_{Z_{\alp,2}}(U_\alp))\subset scl^*(\cf_{Z_2}(X))$.
Hence $\psi\in scl^*(\cf_{Z_2}(X))$. \qed

\hfill

{\bf Proof of Lemma \ref{L:002}.} (1) This was proved in \cite{alesker-mink-oper} in a more general form.

\hfill

(2) We have
\begin{eqnarray}\label{E:jul-1}
F_{\sigma}(\lam,\eps)= \int_{y\in \lam K+\eps A} \sigma([\lam K\cap(y-\eps A)]+B) d vol(y).
\end{eqnarray}

Obviously there exists a constant $C$ such that for any $\lam,\eps\in [0,1]$
$$|\sigma([\lam K\cap (y-\eps A)]+B)|\leq C.$$
Hence for $\lam,\eps\in [0,1]$ one has
$$|F_\sigma(\lam,\eps)|\leq C vol(\lam K+\eps A)=\sum_{j=0}^nC_j \eps^j\lam^{n-j},$$
where $C_j$ are some constants. This implies that the Taylor expansion of $F_\sigma$ at $(0,0)$ does not contain
monomials $\eps^a\lam^b$ with $a+b<n$. This implies part (2) of the lemma.

\hfill

(3) It was shown (in a more general form) in \cite{alesker-mink-oper} that if a sequence $\{\sigma_N\}\subset C^\infty$
converges to $\sigma$ in $C^\infty$ (i.e. uniformly on compact subsets of $V$ with all derivatives)
then
$$F_{\sigma_N}\to F_\sigma \mbox{ in } C^\infty([0,\infty)^2) \mbox{ as } N\to \infty.$$
Hence to prove part (3) of the lemma it suffices to assume that $\sigma$ has a polynomial density on $V$.
We may and will assume that
$$\sigma = P d vol, $$
where $P$ is a {\itshape homogeneous} polynomial of certain degree $d$. Define the function
$$\Phi(\lam,\eps,\delta):=\sigma(\Delta(\lam K)+(\delta B\times \eps A)), \, \lam,\eps\geq 0.$$
By \cite{khovanskii-pukhlikov} (see also \cite{alesker-mink-oper}) this function $\Phi$ is a polynomial
in $\lam,\eps,\delta\geq 0$. Obviously it is homogeneous of degree $d+2n$. Let us write it
$$\Phi(\lam,\eps,\delta)=\sum_{p,q,r}\Phi_{pqr}\lam^p\eps^q\delta^r,$$
where $p,q,r$ must satisfy
\begin{eqnarray}\label{E:jul-2}
p+q+r=d+2n.
\end{eqnarray}

Furthermore $F_\sigma(\lam,\eps)=\Phi(\lam,\eps,1)$ is a polynomial, hence let us write it
$$F_\sigma=\sum_{p,q}F_{pq}\lam^p\eps^q.$$
For the quantity we have to compute we clearly have
\begin{eqnarray}\label{E:jul-5}
\lim_{\lam\to +0}\frac{h_{\sigma,i}(\lam)}{\lam^i}={n \choose i}^{-1}F_{i,n-i}.
\end{eqnarray}
The identity $F_\sigma(\lam,\eps)=\Phi(\lam,\eps,1)$ immediately implies
\begin{eqnarray}\label{E:jul-6}
F_{i,n-i}=\Phi_{i,n-i,d+n}.
\end{eqnarray}

To compute the last expression, let us write
\begin{eqnarray*}
\Phi(\lam,\eps,\delta)=\int_{y\in \lam K+\eps A} \sigma([\lam K\cap (y-\eps A)]+\delta B)d vol(y)=\\
\int_{y\in \lam K+\eps A} \left(\delta^{d+n}\sigma(B)+ (\mbox{ lower degree terms in } \delta)\right) d vol(y)=\\
\sigma(B)\cdot vol(\lam K+\eps A)\cdot \delta^{d+n}+ (\mbox{ lower degree terms in } \delta).
\end{eqnarray*}
This immediately implies that
\begin{eqnarray}\label{E:jul-7}
\Phi_{i,n-i,d+n}={n\choose i}\sigma(B)V(K[i],A[n-i]).
\end{eqnarray}
Lemma follows from (\ref{E:jul-5})-(\ref{E:jul-7}). \qed


\hfill

Let $V_c^{-\infty}(X)$ denote the space of compactly supported generalized valuations. In Part IV we have defined the topology
on this space as follows. For a closed subset $S\subset X$ denote by $V^{-\infty}_S(X)$ the subspace of generalized valuations
with support in $S$. Let us equip $V^{-\infty}_S(X)$ with the weak topology induced from $V^{-\infty}(X)$. Then clearly
$$V_c^{-\infty}(X) =\underset{S\, compact}{\underset{\longrightarrow }{\lim}} V^{-\infty}_S(X)$$
as a vector space. Let us equip $V_c^{-\infty}(X)$ with the topology of inductive limit. Clearly this inductive limit is strict,
namely if $S_1\subset S_2$ then $V_{S_1}^{-\infty}(X)$ is closed in $V_{S_2}^{-\infty}(X)$. Moreover this limit can be taken to be countable
if we take it over a sequence of compact subsets exhausting $X$. Let us denote this topology of inductive limit on $V^{-\infty}_c(X)$ by $\ome$.

Next we have a perfect paring
\begin{eqnarray}\label{E:ind1}
V_c^{-\infty}(X)\times V^\infty(X)\to \CC
\end{eqnarray}
given by $(\phi,\psi)\mapsto \int_X\phi\cdot\psi$.
It is easy to see that this bilinear map is separately continuous.
\begin{claim}
The induced map
$$V_c^{-\infty}(X)\to (V^\infty(X))^*$$
is an isomorphism of vector spaces.
\end{claim}
Under this identification, let us denote by $\sigma$ the weak topology on $V_c^{-\infty}(X)\simeq (V^\infty(X))^*$.
\begin{proposition}
The topologies $\ome$ and $\sigma$ on $V_c^{-\infty}(X)$ coincide.
\end{proposition}
{\bf Proof.} Since the paring (\ref{E:ind1}) is separately continuous when $V^{-\infty}_c(X)$ is equipped with topology $\ome$, it follows
that $\ome$ is finer than $\sigma$. Now let us show that $\sigma$ is finer that $\ome$. We have to show that the identity map
$$(V_c^{-\infty}(X),\ome)\to (V_c^{-\infty}(X),\sigma)$$
is continuous. Equivalently, one has to show  that for any compact subset $S\subset X$ the inclusion map
$$(V_S^{-\infty}(X),\sigma)\to (V_c^{-\infty}(X),\sigma)$$
is continuous. But this is obvious. \qed

\end{document}